\documentclass[a4paper,12pt]{article}

\usepackage[T2A]{fontenc}
\usepackage[cp1251]{inputenc}
\usepackage[russian,english]{babel}

\usepackage{amsfonts}

\def\ds{\displaystyle}
\def\pd{\partial}

\def\a{\alpha}
\def\b{\beta}
\def\d{\delta}
\def\s{\sigma}
\def\t{\tau}
\def\tp{\tau_p}
\def\r{\rho}

\def\E{{\bf E}}
\def\Ex{\E^x}

\def\F{{\cal F}}
\def\R{{\mathbb R}}
\def\L{{\mathbb L}}

\def\G{{\mathcal G}}
\def\M{{\mathcal M}}
\def\I{{\mathcal I}}
\def\II{\mbox{\bf 1}}

\def\beq{\begin{equation}}
\def\eeq{\end{equation}}

\title{\bf Variational View to Optimal Stopping Problems for Diffusion Processes and Threshold Strategies\thanks {The work was supported by Russian
Foundation for Basic Researches (projects 15-06-03723, 13-01-00784) and Russian Foundation for Humanities (project 14-02-00036).
The authors are thankful to Ernst Presman for helpful remarks and discussions.}
}
\author{V.I. Arkin\thanks{Central Economics and Mathematics Institute,
Russian Academy of Sciences.
117418, Moscow,  Nakhimovskii prospect, 47. E-mail:
arkin@cemi.rssi.ru}, A.D. Slastnikov\thanks{Central Economics and Mathematics Institute,
Russian Academy of Sciences.
117418, Moscow,  Nakhimovskii prospect, 47. E-mail: slast@cemi.rssi.ru} }
\date{}

\begin{document}

\maketitle

\begin{quote}
\noindent{\bf Abstract:}
\hangindent=1cm \hangafter=0  We describe a
variational approach to solving optimal stopping problems
for diffusion processes, as an alternative to the traditional
approach based on the solution of the free-boundary
problem. We study smooth pasting conditions from a variational
point of view, and give some examples when the solution to
free-boundary problem is not the solution to optimal stopping problem. A special attention is paid to threshold strategies which allow reduce optimal stopping problem to more simple one-parametric optimization. Necessary and sufficient conditions for threshold structure of optimal stopping time are derived. We apply these results to both investment timing and optimal abandon models.
\end{quote}

\section{Introduction}

Let $X_t,\ t\ge
0$ be a diffusion process with values in $D\subseteq{\mathbb R}^n$ defined
on a stochastic basis $(\Omega, \F, \{\F_t, t \ge 0\},{\bf P})$.

Let us consider an optimal stopping problem for this process:
\beq
\label{optstop}
U(x)=\sup_{\t\in\M} \Ex g(X_\t) e^{-\r\t}\II_{\{\t{<}\infty\}}  ,
\eeq
where
$g:\ D \to \R^1$ is payoff function, $\r> 0$ is discount rate,
and $\Ex$ means the expectation for the process $X_t$ starting from the
initial state $x$. The maximum in (\ref{optstop}) takes over some class $\M$ of stopping times (s.t.)\footnote{In this paper we consider stopping times which can take infinite values (with positive probability)} $\t$. Usually (in classic case) $\M$ is the class of all s.t. with respect to the natural filtration $\F_t^X=\s \{X_s,\ 0\le s\le t\},\ t\ge 0$).

The above problem is one of the classical problems in stochastic control theory. It has a long and wide history in literature and applications.

There are two main approaches to solving an optimal stopping problem for diffusion processes.

The first one, usually referred as Markovian (or `mass'), embeds underlying optimal stopping problem into the family of problems  (\ref{optstop}) with \emph{all} possible initial states $x$ of the process $X_t$. In this case to solve problem (\ref{optstop}) means to find the value function $U(x)$ as well as stopping time $\t^*(x)$, at which $\sup$ in
(\ref{optstop}) is attained  (see, e.g., \cite{ShP, O}).

It follows from general theory that $U(x)$ is the smallest excessive (more precisely, $\r$-excessive) majorant of payoff function  $g(x)$ (Dynkin's characterization, see \cite{Dyn}). Developing this approach in one-dimensional case, Dayanik and Karatzas \cite{DK} gave a characterization of the value
function of the optimal stopping problem (\ref{optstop}) as the smallest nonnegative majorant of payoff $g(x)$, which is concave regarded the certain function which is specified by local characteristics of a given diffusion process. Presman \cite{Pr} proposed how to derive value function using sequential modification of payoff function.

Knowing the value function $U(x)$ one can find the optimal stopping time $\t^*(x)$ as the first exit time of process $X_t$ out of the `continuation set' $C=\{x:\, U(x)>g(x)\}$.

The traditional method for finding the value function for optimal stopping problem (\ref{optstop}) is based on solving the related free-boundary problem:  to find  unknown function
$H(x),\ x\in D$ and set $C\subset \R^n$ such that
\begin{eqnarray}
&&\L H(x)=\r H(x), \quad x\in C;\label{Stefan1}\\
&& H(x)= g(x), \quad x\in \pd C;\label{Stefan2}\\
&& {\rm grad\,} H(x)={\rm grad\,} g(x), \quad x\in \pd
C\label{Stefan3}
\end{eqnarray}
where $\L$ is the infinitesimal generator
 of $X_t$, $\pd C$ is the boundary of the set $C$. The condition
(\ref{Stefan2}) is called ``continuous pasting'', and
 (\ref{Stefan3}) --  ``smooth pasting'' condition.

The solution to the above free-boundary problem is considered as a candidate for solution to optimal stopping problem. And then, using some verification arguments it is necessary to prove that solution $(H(x),\, C)$ to the free-boundary problem really provides the solution of the optimal stopping problem. Namely, $H(x)$ will be a value function and $\t^*(x)=\inf\{t\ge 0: X_t\notin C\}$ will be an optimal stopping time. This free-boundary approach is described in a lot of papers and textbooks, for example, in the monograph by Peskir and Shiryaev \cite{ShP}.

Another approach deals with solving an optimal stopping problem (\ref{optstop}) for fixed initial state $X_0=x$. And, first of all, it should be noted the martingale approach by Beibel and Lerche \cite{BL97, BL01}. Their basic idea is to represent the functional in the right-hand of optimal stopping problem (\ref{optstop}) as a product of positive martingale and `gain' function. And then, using martingale arguments it can be shown that maximization of the gain function gives (under enough weak assumptions) a solution to initial optimal stopping problem. Other applications of martingale methods to solve optimal stopping problems one can find, e.g., in  \cite[Ch. 1]{ShP}.

The present paper develops so-called variational approach to solving an optimal stopping problem for fixed initial state of the process $X_t$, described in \cite{AS09}. In the framework of this approach we propose to find a solution to the problem (\ref{optstop}) over the class $\M$ of stopping times, which are the first exit times of the process $X_t$ from the sets belonging to the given family, and to make optimization over this family of sets.  As an argument for such a reduction of the class of stopping times may be the fact that under enough general assumptions an optimal stopping time in problem (\ref{optstop}) can be find as the first exit time of the process $X_t$ out of the open set $C=\{U(x)>g(x)\}$. Hence we can take $\M$ as a family of first exit times from all open sets (in ${\mathbb R}^n$). Since in one-dimensional case any open set can be represented as countable union of disjoint intervals, then  optimal stopping problem can be reduced (in this case) to finding optimal first exit time from intervals $(a,b),\ l\le a< b\le r$, which contains starting point $x$ of the process $X_t$. Necessary conditions for optimality of such an interval were obtained in \cite{Al01}. The above mentioned intervals are an example of the fact that if the family of sets is chosen `well', then maximization in first exit times from these sets can give a solution to optimal stopping problem over all stopping times. And conversely, one can formulate `an inverse problem':  Under fixed family of sets to describe all conditions on the process and the payoff function, under which an optimal solution (over all stopping times) is the first exit time of the process from a set belonging to the given family.

Note, that for a lot of optimal stopping problem, especially in multi-dimensional case, it is impossible to derive an explicit solution, while the exact solution is very hard for calculations and not suitable for further analysis. Thus, if finding an optimal stopping decision is not the final goal of study (for example, in investment models, see \cite{AS07}), then it makes sense to restrict considerations to simple class of stopping times in order to obtain any `reasonable' solution which will be tractable and suitable for analysis.

At last, the variational approach gives a different look (compared with traditional one) to smooth-pasting principle (\ref{Stefan3}). In the framework of this approach a
smooth-pasting condition can be viewed as first-order optimality
condition for a certain function, while an optimal stopping strategy is associated with a maximization of that function. Thus, a difference between  stationary point and point of maximum can give some non-exotic examples when there are many solutions to free-boundary problem and solution
to free-boundary problem does not give a solution to optimal stopping problem.

The paper is organized as follows. In Section 2  we describe a variational approach for solving an optimal stopping problem.
For the case, when the underlying class of optimal stopping is first exit times from one-parametric family of sets in $\R^n$, an optimal stopping problem can be reduced to one-dimensional maximization of some function. For this case we give necessary and sufficient conditions for optimality of stopping time over the class under consideration.

Frequently, it is optimal to stop when the process exceeds some level (threshold strategy). Similar threshold decisions arise, e.g., in mathematical finance \cite{Shbook}, investment models under uncertainty (real option theory) \cite{DP}, etc. Almost all known decisions in real options theory have a threshold structure (see \cite{DP}). For example, solutions `to invest or not', `to abandon or not' depend on whether the observed values (which determine a decision) will be more or less than some level (threshold).
In Section 3 it is demonstrated how a variational approach works for one-dimensional diffusion processes and two classes of one-parametric sets ($l$- and $r$-intervals), which are generated by threshold strategies.
We give here necessary and sufficient conditions under which the optimal stopping time will have a threshold structure. At the end of this section we apply the obtained results to two fundamental models in real options theory: investment timing model and optimal abandonment model.

In Section 4 we return to general variational approach for solving an optimal stopping problem and consider free-boundary problem (for threshold case) from a variational point of view. Namely, we give an example when a solution to free-boundary problem is not  a solution to optimal stopping problem. Moreover, using second-order optimality conditions we prove some results about a relation between solutions to free-boundary problem and optimal stopping problem.

\section{Variational approach to optimal stopping\\ problem}

In this section we develop approach to solving an optimal
stopping problem which we shall refer as a variational (see \cite{AS09}).
In the framework of this approach, one can define a priori a class $\M$ of stopping times  which are the first exit time out of the set (from a given family of sets), and find optimal stopping time over this class. Besides, unlike the mass setting of an optimal stopping problem, we study the individual problem (\ref{optstop}) for the given (fixed)
initial state of the process $X_0=x$.

Let $\G=\{G\}$ be a given family of regions in $\R^n$,
$\t_G=\t_G(x)=\inf\{t\ge 0:\ X_t\notin G\}$ be a first exit time
of process $X_t$ out of the region $G$ (obviously, $\t_G=0$
whenever $x\notin G$), and $\M({\G})=\{\t_G,\ G\in{\G}\}$ be a set
of first exit times for all regions from the class $\G$.

Under fixed initial value $x$ for any  region
$G\in {\cal G}$ we define the following function (of sets)
\beq
\label{VG}
V_G(x)=\Ex
g(X_{\t_G})e^{-\r\t_G} \II_{\{\t_G{<}\infty\}}.
\eeq
Outside the region $G$ this function equals payoff function $g$, and inside $G$ the function $V_G(x)$ can be derived (under some weak
assumptions) as a solution to boundary Dirichlet problem (see, for example, \cite{Dyn},  \cite{O}):
\beq
\label{Dir}
\begin{array}{ll}
\L u(x)=\r u(x), &x\in
G;\\
 u(y)\to g(x), & y\in G,\ y\to x\in \pd G.
\end{array}
\eeq

In order to calculate functions of the type (\ref{VG}) one can use
martingale methods also (see, for example, \cite{ShP},
 \cite{GSh}).

Thus, a solving an optimal stopping problem (\ref{optstop}) over a class of s.t.
$\M=\M({\cal G})$ can be converted to a solving the
following \emph{variational} problem:
\beq
\label{var}
V_G(x) \to \sup_{G\in {\cal G}}.
\eeq

If $G^*$ is an optimal region in (\ref{var}), then the  first exit time from this region  $\t_{G^*}$ will be the optimal stopping time for the problem (\ref{optstop}) over the class $\M=\M(\G)$.

If the class of regions $\cal G$ is chosen `well', it is possible to prove that s.t. $\t_{G^*}$ will  be  also an optimal  stopping
time for problem  (\ref{optstop}) over all s.t.
 $\M$. In \cite{AS09} such approach was realized  for  two-dimension
geometric  Brownian  motion  $X_t$ and  homogeneous payoff function
$g$.

A close approach is developed in \cite{Al01}, where an optimal
stopping problem for one-dimensional diffusion is solved by
mathematical programming technique. However those method uses a
few properties of one-dimensional diffusion, which are not valid
in multi-dimensional case.

As for multi-dimensional processes, the first-order conditions as
a heuristic method for finding boundaries of optimal
``continuation sector'' in optimal stopping problem for bivariate
geometric Brownian motion and homogeneous (of first order) payoff
function was used in \cite{GSh}.

Let  us  note, that the calculation of the optimal stopping time over
a given class  of regions represents, to our opinion, a practical
interest. Indeed, free-boundary problem for multi-dimensional diffusion processes has no (as a rule) explicit solution. Therefore, it has
a sense to restrict our consideration to more simple stopping times and corresponding regions, for
which it is possible to derive the function of sets $V_G(x)$. Also,
numerical methods can be applied for solving the problem (\ref{var})
with fixed initial state $X_0=x$.

An idea of variational approach is general enough and can be applied
not only for a diffusion processes and payoff functions of the type
(\ref{optstop}).

\subsection{One-parametric family  of regions}

Under some additional assumptions a general variational problem
(\ref{var}) can be simplified and be made more convenient for study.

Let $\G=\{G_p,\ p\in
P\subset \R^1\}$ be one-parametric family of regions in $\R^n$,
$\t_p=\inf\{t{\ge} 0: X_t\notin G_p\}$,
\beq
\label{Vpx}
V(p;x)=V_{G_p}(x)=\Ex
g(X_{\t_p})e^{-\r\t_p} \II_{\{\t_p{<}\infty\}}
\eeq
(see formula (\ref{VG})).

The function $V(p;x)$ is defined on $P\times D$,  and, obviously,
$V(p;x)=g(x)$ for $x\notin G_p$.

Further, we assume that a family of regions $\{G_p\}$ satisfies the
following conditions:
\begin{description}
\item[(A1)] \emph{Monotonicity}: \ $G_{p_1}\subset G_{p_2}$ whenever $p_1<p_2$.
\item[(A2)] \emph{Thickness}: \  Every point $x\in D$ belongs to the boundary of the unique set from
the class $\G$. The parameter of those set will be refered as $q(x)$, so $x\in \pd G_{q(x)}$.
\end{description}

Under a thickness property one can write:
\beq
\label{cont}
V(q(x);x)=g(x) \quad \forall x\in D.
\eeq

For a one-parametric case the variational problem (\ref{var}) can be rewritten as one-dimensional optimization:
\beq
\label{var1}
V(p;x) \to \sup_{p\in {P}},
\eeq
where $V(p;x)$ is specified in (\ref{Vpx}).

Under the assumptions \textbf{(A1)} and \textbf{(A2)} necessary and sufficient conditions for  maximization of $V(p;x)$ in $p$ are given by the following
\medskip

{\bf Theorem 1.} {\it

\emph{i)}
If $p^*=p^*(x)$ is the solution to the problem (\ref{var1}) then conditions
\begin{eqnarray}
&&
V(p;x)\le V(p^*;x)\  \mbox{\rm whenever }  p<p^*,\ x\in G_p\cup\pd G_p, \label{criteria-var0}\\
&&V(p;x)\le g(x) \ \mbox{\rm whenever } p>p^*,\ x\in G_p\setminus G_{p^*}\label{criteria-var00}
\end{eqnarray}
hold.

\emph{ii)}
If for some $p^*=p^*(x)$
\beq\label{criteria-var1}
V(p_1;x)\ge V(p_2;x)\ \mbox{\rm whenever }  p^*\le p_1 <p_2,\ x\in G_{p_2}
\eeq
and condition (\ref{criteria-var0}) hold, then $p^*$ is the solution to the problem (\ref{var1}).
 }
\medskip

{\sc Proof. } \
i) Inequality in (\ref{criteria-var0}) is the direct consequence of optimality of $p^*$ in problem (\ref{var1}). For $p>p^*$ and $x\in G_p\setminus G_{p^*}$ we have $g(x)=V(p^*;x)\ge V(p;x)$, i.e. (\ref{criteria-var00}).

ii) Show that
$V(p;x)\le V({p^*};x)$ for any $p\in P$. Denote $\overline{G}_p=G_p\cup\pd G_p$ --- the closure of the set $G_p$.

 Let $p<p^*$. If $x\notin G_{p^*}$, then $x\notin G_{p}$ (due to
monotonicity of regions) and, therefore, $V(p;x)= g(x)=
V({p^*};x)$. If $x\in G_{p^*}$ and $x\in \overline{G}_p$, then $V(p;x)\le V(p^*;x)$ (due to (\ref{criteria-var0})). Finally, if $x\in G_{p^*}$ and
$x\notin \overline{G}_p$, then $q(x)<p^*$, $x\in \pd G_{q(x)}$, and  using~(\ref{cont}), (\ref{criteria-var0}), we have: $V(p;x)=
g(x)=V(q(x);x)\le V(p^*;x)$.

 Consider the case $p>p^*$. If $x\notin G_{p}$, then $x\notin
G_{p^*}$ (due to monotonicity of regions), hence, $V(p;x)= g(x)=
V({p^*};x)$. Whenever $x\in G_{p}$ and ${x\in G_{p^*}}$, then
$V(p;x)\le V(p^*;x)$ due to (\ref{criteria-var1}). When $x\in G_{p}$ and
$x\notin G_{p^*}$ one can see that $p^*\le q(x)< p$. Therefore, (\ref{criteria-var1}) implies: $V(p;x)\le V(q(x);x)=g(x)=
V({p^*};x)$. $\square$
\medskip

As one can see there is a gap between necessary conditions (\ref{criteria-var0}), (\ref{criteria-var00}) and sufficient conditions (\ref{criteria-var0}), (\ref{criteria-var1}). But for one-dimensional diffusion processes  it will be derived that the their direct consequences  give necessary and sufficient conditions for optimality of stopping time, so the mentioned above gap disappears (see Theorem 2 in Section 3).

\subsection{A variational look to smooth pasting principle}

The variational approach can give a new look to a smooth pasting principle, which is crucial in solving free-boundary problem.

Let the set of states $D\subseteq \R^n$ of the diffusion process $X_t$ be an open set, the assumptions (A1)--(A2) on one-parametric family of sets $\G=\{G_p,\, p\in P\}$, where $P$ is an open set in $\R^1$, hold. Suppose that functions $g(x)$, as well as $q(x)$ (the parameter of the region whose boundary passes through the point $x$), are
differentiable in all arguments. Note, that the function $q(x)$
will be smooth, for example, in the case when regions' boundaries
 $\pd G_p$ are specified by surfaces of the type $\{\Psi(p,y)=0,\
 y\in \R^{n}\}$, where $\Psi(p,y)$ is continuously differentiable
 in $(p,y)$ and  $\ds \Psi'_p(p,y)$ is non-zero.
Moreover, assume there exists differentiable function $F(p,x)$ on $P\times D$ such that $F(p,x)=V(p;x)$ for $p\in P,\ x\in G_p$ (where $V(p;x)$ is defined in (\ref{Vpx})).

Let $\bar p(x)$ be a stationary point of the function $F(p,x)$ in
$p$, i.e. \linebreak $F'_p(\bar p(x),x) = 0$ ($x\in D$). The continuous pasting condition $V(q(x);x)=g(x)$ (see (\ref{cont})) implies
\beq
\label{smooth0}
F_p'(q(x),x)\
{\rm grad}\, q(x) + F'_x(q(x),x) = {\rm grad\,} g(x), \quad x\in D.
\eeq
Thus, if $x\in \pd G_{\bar p(x)}$, then $q(x)=\bar p(x)$ and, therefore,
\beq
\label{smooth1}
F'_x(\bar
p(x),x) = {\rm grad}\,g(x), \quad x\in \pd G_{\bar p(x)}.
\eeq

This equality can be viewed as a variant of smooth pasting condition at the
boundary of the set $G_{\bar p(x)}$, whose parameter is a stationary
point of a function  $F(p,x)$.

Note, the set of such $x$ that (\ref{smooth1}) holds can be empty.
Consider the case, when stationary points $\bar p(x)=\bar p$ do not
depend on $x$. In this case the set of such $x$ that relation
(\ref{smooth1}) valid, is not empty. As we see below at Section 3, such a situation emerges, in
particular, for one-dimensional diffusions and the classes of $l$-intervals or $r$-intervals. In these cases the function $F(p,x)$  has a multiplicative structure: $F(p,x)=F_1(p)F_2(x)$ --- see (\ref{representation}), (\ref{two-fpx}).

Defining the function $\overline F(x)=F(\bar p,x)$, the relation
(\ref{smooth1}) can be written as follows:
\beq
\label{smooth2}
{\rm
grad}\, V_{\bar p}(x)={\rm grad}\,\overline F(x) = {\rm grad}\,g(x), \quad
x\in \pd G_{\bar p}.
\eeq

Taking into account that $\overline F(x)$ for $x\in G_{\bar p}$ is a
solution to Dirichlet problem (\ref{Dir}), the equality
(\ref{smooth2}) is a traditional smooth pasting condition, and,
therefore, the pair $(\overline F(x),G_{\bar p})$ is a solution to free-boundary problem (\ref{Stefan1})--(\ref{Stefan3}).
Let us note, if ${\rm grad}\,q(x)\neq 0$ for the family of regions $\{G_p\}$ and $x\in D$, then (as one can see from
(\ref{smooth0})), the smooth pasting condition (\ref{smooth2}) is
equivalent to stationarity of $F(p,x)$ in $p$ at the point
$\bar p$.

On the other hand, if maximum of the
function $F(p,x)$ in $p$ is attained at the point $p^*\in P$, then the first exit time $\t_{p^*}$ from the region $G_{p^*}$ will
be a candidate for an optimal stopping time over the class $\M(\G)$ (see Theorem 1).

Of course, point of maximum is a stationary point of the function $F$, but not vice versa. Hence, for such a case a solution to free-boundary problem can not give a solution to optimal stopping problem. We continue to discuss this question in Section 4.

\section{One-dimensional diffusion processes.\\ Threshold strategies}

Let $X_t$ be a diffusion process with values in the  segment $D\subseteq \R^1$ with boundary points $l$ and $r$, where $ -\infty \le l< r\le +\infty$, open or closed (i.e. it may be $(l,r)$,\, $[l,r)$,\, $(l,r]$,\, or $[l,r]$),  and its
infinitesimal generator has the following type :
\beq \label{infop1}
\L f(x)=a(x)f'(x)+ {\textstyle\frac12} \, \s^2(x) f''(x),
\eeq
where $a:\ D \to \R^1$ and  $\s:\ D \to \R^1_+$ are the drift and diffusion functions, respectively. Denote $\I=(l,r)$.

The process $X_t$ is assumed to be regular; this means
that, starting from an arbitrary point $x\in \I$, this process
reaches any point $y\in \I$ in finite time with positive
probability.
To guarantee the regularity of $X_t$ we assume that the drift and diffusion functions are satisfy the following local integrability condition:
\beq \label{regular}
 \int_{x-\varepsilon}^{x+\varepsilon}\frac{1+|a(y)|}{\s^2(y)} dy <\infty \quad \mbox{for some }\varepsilon>0,
\eeq
at any $x\in \I$ (see, e.g. \cite{KS}).
\smallskip

It is known that under regularity conditions (\ref{regular}),
on the interval $\I$, there exist (unique up to constant positive multipliers) increasing and decreasing functions $\psi(x)$ and $\varphi(x)$ with absolutely continuous derivatives, which are the fundamental  solutions to the ODE
\beq\label{diffur}
\L u(x)=\r u(x)
\eeq
almost sure (in Lebesque measure) on the interval $\I$ (see, e.g. \cite[Chapter 5, Lemma 5.26]{KS}). Moreover, $0<\psi(x),\, \varphi(x)<\infty$ for $x\in \I$.
Note, if functions $a(x),\ \s(x)$ are continuous, then \ $\psi,\, \varphi\in C^2(\I)$.

Below we introduce two natural one-parametric family of sets  underlying a variational approach for one-dimensional diffusions.

\subsection{$l$-intervals}

As the first one-parametric family of sets we take intervals of the type
$G_p=\{x\in D:\, x<p\}$, \ $p\in (l,r)$, which we call $l$-intervals (emphasizing that \emph{left} end is fixed and equal to left boundary point of the process values $D$).  $l$-interval $G_p$ is $[l,p)$ or $(l,p)$ in dependence on $l\in D$ or not. Obviously, the class of $l$-intervals satisfies conditions {\rm{\bf (A1)}--{\bf (A2)}} in Section 2.1.

Let us define $\tp^l=\inf\{t{\ge} 0: X_t\notin G_p\}=\inf\{t{\ge} 0: X_t\ge p\}$ --- the first time when the process $X_t$ leaves  $G_p$.  We will call $\tp^l$ as threshold stopping time (first exit time \emph{over} threshold $p$). Let $\M^l_{\rm th}=\{\tp^l, \ p\in \I\}$ be a class of all such threshold stopping times.

Along with the above stopping time let us define the first hitting time to threshold:  $T_p= \inf\{t{\ge} 0: X_t= p\},\ p\in (l,r)$.

Then for the above class of $l$-intervals the function $V(p;x)$, defined in previous section, has the following representation:
\medskip

{\bf Lemma 1.} {\it If  $x,p\in \I$, then
\beq\label{representation}
V(p;x)= \left\{
\begin{array}{ll}
g(p)\psi(x)/\psi(p), & \mbox{\rm for } x<p,\\[4pt]
g(x), &    \mbox{\rm for } x\ge p,
\end{array}
\right.
\eeq
where $\psi(x)$ is an increasing solution to ODE (\ref{diffur}).
}
\medskip

{\sc Proof. } \
Due to known formulas (see, e.g., \cite{IM, BorSal}):
$$
\Ex e^{-\r T_p}= \left\{
\begin{array}{ll}
\psi(x)/\psi(p), & \mbox{\rm for } x<p,\\[4pt]
\varphi(x)/\varphi(p), &    \mbox{\rm for } x\ge p.
\end{array}
\right.
$$

Therefore, for $x<p$, obviously, $\tp^l=T_p$ and using the above  formula we have:
$$
\Ex g(X_{\tp^l})e^{-\r\tp^l} \II_{\{\tp^l{<}\infty\}}=\Ex g(X_{T_p})e^{-\r T_p}=g(p)\Ex e^{-\r T_p} =g(p)\frac{\psi(x)}{\psi(p)},
$$
that proves lemma. $\square$
\medskip

For one-dimensional case we can enhance Theorem 1 and give necessary and sufficient conditions (which are the same) for optimality over the class of $l$-intervals (or threshold stopping times). Let us define
\beq \label{h}
h(p)= g(p)/\psi(p).
\eeq

In order to exclude trivial cases (when optimal threshold stopping time does not exist) we assume that $g(x_0)>0$ for some $x_0\in \I$.
\medskip

{\bf Theorem 2.} {\it
Threshold stopping time $\t^l_{p^*}$ is optimal in the problem (\ref{optstop})  over the class $\M{=}\M^l_{\rm th}$ for all $x\in \I$, if and only if the following conditions hold:
\begin{eqnarray}
&& h(p)\le h(p^*)\  \mbox{\rm whenever }  p<p^*; \label{criteria0}\\
&& h(p)\ \mbox{\rm does not increase for }  p\ge p^* \label{criteria00}.
\end{eqnarray}
}
%\smallskip

{\sc Proof. } \ Using Lemma 1 it is easy to see that in our case the condition  (\ref{criteria-var0}) is equivalent to (\ref{criteria0}), and (\ref{criteria-var00}) is equivalent to (\ref{criteria00}). Therefore, by Theorem 1, optimality of  $\t^l_{p^*}$ implies conditions (\ref{criteria0})--(\ref{criteria00}).

On the other hand, let (\ref{criteria0})--(\ref{criteria00}) hold. Take arbitrary $p_1,p_2$ such that $p^*\le p_1<p_2$ and $x<p_2$. If $x<p_1$ then $V(p_1;x)=h(p_1)\psi(x)\ge h(p_2)\psi(x)=V(p_2;x)$. And if $p_1\le x<p_2$ then $V(p_1;x)=g(x)=h(x)\psi(x)\ge h(p_2)\psi(x)=V(p_2;x)$. Thus, the condition  (\ref{criteria-var1}) holds and we can apply Theorem 1 again. $\square$

\medskip

So, the optimal threshold $p^*$ is a point of maximum for the function $h(p)$. This implies the necessity (under minor assumptions) of the `generalized' smooth-pasting principle.
\medskip

{\bf Corollary 1.} {\it Let $\t_{p^*}$, $p^*\in \I$, be the optimal stopping time in the problem (\ref{optstop}) over the class $\M^l_{\rm th}$, and function  $g(x)$ has one-sided derivatives $g'(p^*{\pm}0)$  at the point $p^*$. Then the function $\ds v(x)=\sup_{p\in \I} V(p;x)$ has one-sided derivatives at the point $p^*$ and the following inequalities hold:
$$
 g'(p^*{+}0)=v'(p^*{+}0)\le v'(p^*{-}0)\le g'(p^*{-}0).
$$
}

{\small\sc Proof. } \ From (\ref{representation}) it follows \beq\label{representation1}
v(x)=V(p^*;x)=\left\{
\begin{array}{ll}
h(p^*)\psi(x), & \mbox{\rm for } x<p^*,\\[4pt]
g(x), &    \mbox{\rm for } x\ge p^*.
\end{array}
\right.
\eeq

Since $p^*$ is a point of maximum for the function $h(p)$, then
$
h'(p^*{-}0)\ge 0\ge h'(p^*{+}0).
$
Therefore,
$$
v'(p^*{-}0)=h(p^*)\psi'(p^*)=g'(p^*{-}0)-h'(p^*{-}0)\psi(p^*)\le g'(p^*-0),
$$
$$
v'(p^*{+}0)= g'(p^*{+}0)=h'(p^*{+}0)\psi(p^*)+h(p^*)\psi'(p^*)\le h(p^*)\psi'(p^*)=v'(p^*{-}0).
$$

The corollary proved.
\medskip

{\bf Remark.} {\it If function $g(x)$ is differentiable at the point $p^*$, then the function $v(x)$ will be differentiable at the
point $p^*$, and ${v}'(p^*){=}g'(p^*)$.
}
\medskip

The necessity of the smooth pasting condition under
some additional constraints on the process was shown
in \cite{ShP}. A result similar to ours was obtained in \cite{V}

The extended conditions of Theorem 2 will be necessary and sufficient for that a solution to stopping problem (\ref{optstop}) over threshold stopping times remains optimal over \emph{all} stopping times.

Let us define for $p\in \I$ a restriction of a process $X_t$ on $D^l_p=\{x\in D:\, x\ge p\}$ as $Y_t=X_{t\wedge T_p}$ where $Y_0=x>p$ and $T_p$ is the first hitting time of threshold $p$. In other words, the process $Y_t$ is started from a point greater than $p$ and is absorbed when it reaches $p$.

We will call the nonnegative function $f:\, D^l_p\to {\mathbb R}_+$ as $\r$-excessive with respect to restriction of $X_t$ on $D^l_p$ if for all $x>p$ and any stopping time $\t$
\beq\label{ex}
\Ex f(X_{\t\wedge T_p}) e^{-\r(\t\wedge T_p)}\II_{\{\t\wedge T_p{<}\infty\}}\le f(x).
\eeq
where  $T_p= \inf\{t{\ge} 0: X_t= p\},\ p\in (l,r)$.
\medskip

{\bf Theorem 3.} {\it
Let $\t^l_{p^*}$, $p^*\in \I$, be an optimal stopping time in the problem (\ref{optstop}) over the class $\M^l_{\rm th}$ of threshold stopping times for all $x\in \I$, moreover,  there exists $g'(p^*{+}0)$ and $g(x)\ge 0$ for $x>p^*$. 
Then threshold stopping time $\t^l_{p^*}$ is the optimal stopping time in problem (\ref{optstop})  over class of all stopping times for all $x\in \I$, i.e.
$$
U(x)=\sup_{\t}\Ex g(X_\t) e^{-\r\t}\II_{\{\t{<}\infty\}}=\left\{
\begin{array}{ll}
h(p^*)\psi(x), & \mbox{\rm for } x<p^*,\\[4pt]
g(x), &    \mbox{\rm for } x\ge p^*,
\end{array}
\right.
$$
if and only if the function $g(x)$ is $\r$-excessive with respect to restriction of process $X_t$ on $D^l_{p^*}$.
}

{\sc Proof of Theorem 3.} \
Let $\t^l_{p^*}$ be optimal in problem (\ref{optstop}) over all stopping times for all $x\in \I$. Then for $x>p^*$ and any stopping time $\t$ we have
$$
\Ex g(X_{\tilde\t}) e^{-\r\tilde\t}\II_{\{\tilde\t{<}\infty\}}\le U(x)=g(x),
$$
where $\tilde\t=\t\wedge T_{p^*}$. It means that $g(x)$ is $\r$-excessive function with respect to restriction of process $X_t$ on $D^l_{p^*}$.

To establish the inverse implication let define the function
$$
\Phi(x)=V(p^*;x)=\left\{
\begin{array}{ll}
h(p^*)\psi(x), & \mbox{\rm for } x<p^*,\\[4pt]
g(x), &    \mbox{\rm for } x\ge p^*.
\end{array}
\right.
$$
and  show that $\ds \Ex \Phi(X_{\t}) e^{-\r\t}\II_{\{\t{<}\infty\}}\le \Phi(x)$ for any stopping time $\t$ and $x\in \I$.

In order to prove it  we use the following criteria of  $\r$-excessive functions from \cite[Theorem 5.1]{LZ}.
\medskip

{\bf Lemma 2.} {\it
A function $F:\, D\to {\mathbb R}_+$ satisfies the inequality
$$
\Ex F(X_{\t}) e^{-\r\t}\II_{\{\t{<}\infty\}}\le F(x)
$$
for all stopping times $\t$ and all initial states $x\in \I$, if and only if the following statements hold:

i) $F(x)$ is the difference of two convex functions for $x\in \I$;

ii) the measure ${\mathcal L} F (dx)={\textstyle\frac12} \, \s^2(x) F''(dx)+a(x)F'(x{-}0)dx -\r F(x)dx$, where $F''(dx)$ means second distributional derivative, is non-positive;

iii) $F$ is lower semi-continuous at absorbing boundary points.
}
\medskip

Apply this lemma to the function $\Phi(x)$.

If $x<p^*$ then $\Phi(x)=h(p^*)\psi(x)$. Obviously, $\psi(x)$ is the difference of two convex functions, because $\psi'(x)$ has a finite variation (at any segment $[l',p^*],\ l'>l$) and, therefore, is the difference of two increasing functions. Thus, condition i) of Lemma 2 is satisfied. Further, ${\mathcal L} \Phi (dx)=h(p^*)[\L\psi(x)-\r \psi(x)]dx=0$ a.s. (by definition of $\psi(x)$), i.e. ii) holds. At last, if $l$ is absorbing state we can define $\psi(l)$ in continuity, so that condition iii) of Lemma 2 is also satisfied.

If $x>p^*$ then $\Phi(x)=g(x)$. Since $g(x)$ is $\r$-excessive with respect to restriction of process $X_t$ on $D^l_{p^*}$, then conditions i)--iii) of Lemma 2 are satisfied. 

Finally, at the point $x=p^*$ we have ${\mathcal L} \Phi (\{p^*\})={\textstyle\frac12} \, \s^2(p^*)[g'(p^*{+}0)-h(p^*)\psi'(p^*)]\le 0$, since $p^*$ is the point of maximum for function $h(p)$ (see Theorem 2) and, therefore $h'(p^*{+}0)\le 0$.

Thus, all conditions i)--iii) of Lemma 2 are satisfied for the function $\Phi(x),\, x\in D$. 
 
Since  $\t^l_{p^*}$ is optimal over the class $\M^l_{\rm th}$, then Theorem 2 implies $h(p^*)\psi(x)\ge h(x)\psi(x)=g(x)$ for $x<p^*$, therefore $\Phi(x)\ge g(x)$ for all $x$. Using this inequality and Lemma 2 we have
$$
\Ex g(X_{\t}) e^{-\r\t}\II_{\{\t{<}\infty\}}\le \Ex \Phi(X_{\t}) e^{-\r\t}\II_{\{\t{<}\infty\}}\le \Phi(x)
$$
for all stopping times $\t$ and $x\in \I$. Hence, $ U(x)\le \Phi(x)$. 
On the other hand, obviously, $U(x)\ge\Phi(x)$.

Therefore, $ U(x)= \Phi(x)=V(p^*;x)$, i.e. $\t_{p^*}$ is the optimal stopping time in problem (\ref{optstop}) over all stopping times for all $x$. 

This completes the proof.

\medskip

As we see the question ``Does optimal threshold stopping time  $\t_{p^*}$ remain optimal over all stopping times?'' is reduced to problem of $\r$-excessiveness of payoff function with respect to the restriction of underlying process on the range $D^l_{p^*}$. There are several conditions and criteria for  $\r$-excessiveness, see, e.g. \cite{DK, LZ}.   

Now we give necessary and sufficient conditions for such a $\r$-excessiveness of payoff function, which may be convenient for applications.
\medskip

{\bf Statement.} 
{\it 
Let for a set of isolated points  $\{a_1,...,a_n,...\}$, where \linebreak $p^*{<}a_1{<}a_2{<}...{<}r$, the function $g'(p)$ be absolutely continuous on the intervals $(p^*, a_1), \ (a_i, a_{i+1})$, $i{\ge} 1$ and there exist one-sided derivatives $g(p^*{+} 0),\, g'(a_i{\pm} 0)$,\, $i{\ge} 1$, such that $\sum\limits_{i\ge 1} \s^2(a_i) |g'(a_i{+}0){-}g'(a_i{-}0)|<\infty$.
Then $g(x)$ is $\r$-excessive function  with respect to restriction of process $X_t$ on $D^l_{p^*}$ if and only if the following conditions hold:

i) $\L g(p)\le \r g(p)$ \ a.s. (in Lebesque measure) for  $p>p^*$;

ii) $g'(a_i{+} 0) - g'(a_i{-} 0)\le 0$ for all $i\ge 1$.
}
\smallskip

(This statement easy follows from Lemma 2.)

\subsection{$r$-intervals }

Now let us consider another one-parametric family of sets underlying the variational approach to solving optimal stopping problems, namely, the class of $r$-intervals of the type $G_p=\{x\in D:\, x>p \},\, p\in \I$. This class is similar to the class of $l$-intervals, but here we fix the right end of the range of the process $D$ instead of its left end (as in $l$-intervals).

In complete analogy with the above considerations we can define threshold stopping time $\t^r_p=\inf\{t{\ge} 0: X_t\le p\}$  as a first time when process $X_t$ falls \emph{below} threshold $p$, \ $p\in \I$, and the corresponding class $\M^r_{\rm th}$.

Then for the class of $r$-intervals the function $V(p;x)$, defined in (\ref{Vpx}), has the following representation (cf. (\ref{representation})):
\beq\label{two-fpx}
V(p;x)= \left\{
\begin{array}{ll}
g(p)\varphi(x)/\varphi(p), & \mbox{\rm for } x>p,\\[4pt]
g(x), &    \mbox{\rm for } x\le p,
\end{array}
\right.
\eeq
where $\varphi (x)$ is a decreasing solution to ODE (\ref{diffur}).

Thus, we can modify all results from the previous section for the class of $r$-intervals. 
\medskip

{\bf Theorem 2$'$.} {\it
Threshold stopping time $\t^r_{p^*}$ is optimal in the problem (\ref{optstop})  over the class $\M^r_{\rm th}$ for all $x\in \I$, if and only if the following conditions hold:
\begin{eqnarray*}
&& {g(p)}/{\varphi(p)}\le {g(p^*)}/{\varphi(p^*)}\  \mbox{\rm whenever }  p>p^*; \\[2pt]
&& {g(p)}/{\varphi(p)}\ \mbox{\rm does not decrease for }  p\le p^* .
\end{eqnarray*}
}
\medskip

We can  define $D^r_p=\{x\in D:\, x\le p\}$ for $p\in \I$ and a restriction of a process $X_t$ on $D^r_p$ as $Y_t=X_{t\wedge T_p}$, where $Y_0=x<p$ and $T_p$ is the first hitting time of threshold $p$. 

In analogue with the case of $l$-intervals we call the function $f:\, D^r_p\to {\mathbb R}_+$  $\r$-excessive with respect to restriction of $X_t$ on $D^r_p$ if for all $x<p$ and any stopping time $\t$ the inequality (\ref{ex}) hold.

Now, the analogue of Theorem 3 is the following.
\medskip

{\bf Theorem 3$'$.} {\it
Let $\t^r_{p^*}$, $p^*\in \I$, be an optimal stopping time in the problem (\ref{optstop}) over the class $\M^r_{\rm th}$ of threshold stopping times for all $x\in \I$, and  there exists $g'(p^*{-}0)$.
Then threshold stopping time $\t^r_{p^*}$ is the optimal stopping time in problem (\ref{optstop})  over class of all stopping times for all $x\in \I$, i.e.
$$
U(x)=\sup_{\t}\Ex g(X_\t) e^{-\r\t}\II_{\{\t{<}\infty\}}=\left\{
\begin{array}{ll}
g(p^*)\varphi(x)/\varphi(p^*), & \mbox{\rm for } x>p^*,\\[4pt]
g(x), &    \mbox{\rm for } x\le p^*,
\end{array}
\right.
$$
if and only if the function $g(x)$ is $\r$-excessive with respect to restriction of process $X_t$ on $D^r_{p^*}$.
}
\medskip

\medskip

Threshold stopping times $\t^r_p$ appears, in particular, in optimal stopping problems with both integral and terminal payoffs:
\beq\label{integ}
\E^x\left(\int_0^{\t}g_1(X_t) e^{-\r t}\,dt + g_0(X_\t) e^{-\r
\t}\right)\to \sup_\t,
\eeq
where $g_0(x),\, g_1(x)$ are given functions.

Indeed, using strict Markov property of diffusion processes one can reduce the problem (\ref{integ}) to optimal stopping problem (\ref{optstop}) with terminal payoff
\beq\label{reduction}
g(x)=g_0(x)-R(x), \quad \mbox{ where } \ R(x)=\E^x \int_0^\infty g_1(X_t) e^{-\r t}\,dt.
\eeq

By Green function (see, e.g., \cite{BorSal,Al01}) one can get the following Green representation for $R(x)$:
\beq\label{Green}
R(x)=B^{-1}\left(\varphi(x)\int_l^{x}\psi(y)g_1(y)H(y)\,dy + \psi(x)\int_x^r \varphi(y)g_1(y)H(y)\,dy\right),
\eeq
where $B=[\psi'(x)\varphi(x)-\psi(x)\varphi'(x)]/S'(x)>0$ is constant, \ $H(y)=2[\s^2(y)S'(y)]^{-1}$, \ $\ds S'(x)=\exp\left\{-\int \frac{2 a(x)}{\s^2(x)}dx\right\}$ is derivative of the scale function for the process $X_t$.
\medskip

If functions $g_1(x)$ and $a(x),\, \s(x)$ are continuous, then Green representation (\ref{Green}) implies that  $R(x)$ will be twice differentiable at interval $\I$ and
\begin{eqnarray}
&&  R'(x)=B^{-1}[\varphi'(x)I_1(x) + \psi'(x)I_2(x)],\\[3pt]
&&R''(x)= B^{-1}[\varphi'(x)I_1(x) + \psi'(x)I_2(x)]-2g_1(x)/\s^2(x);\\[3pt]
&& I_1(x)=\int_l^{x}\psi(y)g_1(y)H(y)\,dy,\quad
 I_2(x)=\int_x^r \varphi(y)g_1(y)H(y)\,dy,\label{I}
\end{eqnarray}
where $B$ is defined in (\ref{Green}).

Hence, Theorems 2$'$, 3$'$ together with the variant of Statement (for $r$-intervals) imply the following result concerning a threshold structure of solution to  optimal stopping problem with both integral and terminal payoffs (\ref{integ}).
\medskip

{\bf Theorem 4.}
{\it
Let for some $p^*\in \I$ functions $a(x),\s(x), g_1(x)$ be continuous, $g_0(x)$ be twice differentiable and $g_0(x)\ge R(x)$ on segment $(l,p^*]$.
Then threshold stopping time $\t^r_{p^*}=\inf\{t{\ge} 0: X_t\le p^*\}$ is the optimal stopping time in problem
(\ref{integ}) over the class of all stopping times (for all $x\in \I$) if and only if the following conditions hold:
\begin{eqnarray}
&&[g_0(p)-R(p)]\varphi(p^*)\le [g_0(p^*) - R(p^*)]\varphi(p)\quad  \mbox{\rm for } p>p^*;\label{criteria101}\\[3pt]
&& I_2(p^*)S'(p^*)= g'_0(p^*)\varphi(p^*) - g_0(p^*)\varphi'(p^*);\label{criteria1011}\\[3pt]
&&  \L g_0(p)-\r g_0(p)\le \L R(p)-\r R(p) \quad \mbox{\rm   for } p<p^*, \label{criteria12}
\end{eqnarray}
where function $I_2(x)$ is specified in (\ref{I}) and  $\varphi (x)$ is a decreasing solution to ODE (\ref{diffur}).
}
\medskip

Results presented in Sections 3.1 and 3.2 extend the corresponding results obtained in \cite{A15, AS12, CM, V}.

\subsection{Application to Real Options}

\paragraph{\hspace*{\parindent} Investment timing problem.}

One of the fundamental problems in real options theory is to derive the optimal time at which an investor (or decision maker) should finance and launch an investment project --- investment timing problem.

Let $I$ be a cost of investment required for beginning a project, and $X_t$ is considered as Present Value (PV) from the project started at time $t$. As usual investment supposed to be instantaneous and irreversible, and the project --- infinitely-lived.

At any time a decision-maker (investor)  can either {\it accept} the project and proceed  with  the investment  or
{\it delay} the decision until he obtains  new information
regarding its environment (prices of the product and resources,
demand etc.). Thus, the main goal of a decision-maker in this
situation is to find, using the available information, an
optimal time for investing the project (investment timing
problem), which maximizes the Net Present Value from the project:
\beq
\label{optinv0}
{\bf E}\left(X_\t - I\right) e^{-\r \t} \to \max_\t,
\eeq
where the maximum is considered over all stopping times $\t{\in}\M$.

The majority of results on this problem which was introduced in \cite{MS}, has a threshold structure: to invest when PV from the project exceeds the certain level (threshold). In the heuristic level it was shown for the cases when PV is modelled by geometric Brownian motion, arithmetic Brownian motion, mean-reverting process and some other (see \cite{DP}). And the general question arises: For what processes of PV from the project an optimal decision to an investment timing problem will have a threshold structure? Some sufficient conditions in this direction was obtained in \cite{Al01}. We can give the necessary and sufficient conditions for optimality of threshold strategy in investment timing problem (\ref{optinv0}) as the consequence of our above results (but not a direct corollary) for the linear payoff function $g(x)=x-I$. We consider the case $I<r$ else the optimal time in (\ref{optinv0}) will be $+\infty$.

\medskip
{\bf Theorem 5.} {\it
Threshold stopping time  $\t^l_{p^*}$, $p^*\in (l,r)$, is optimal in the investment timing problem (\ref{optinv0}) if and only if the following conditions hold:
\begin{eqnarray}
&&
(p-I){\psi(p^*)}\le (p^*-I){\psi(p)}\quad \mbox{\rm  for } p<p^*;\label{criteria110}\\[3pt]
&& \psi(p^*)=(p^*-I)\psi'(p^*);
\label{criteria11}\\
&& a(p)\le \r (p- I) \quad \mbox{\rm  for } p>p^*, \label{criteria21}
\end{eqnarray}
where $\psi(x)$ is an increasing solution to ODE (\ref{diffur}) and $a(p)$ is the drift function of the process $X_t$.
}
\medskip

{\sc Proof.} \
Let $\t^l_{p^*}$, $p^*\in (l,r)$, is optimal in the investment timing problem. Then it is optimal in (\ref{optinv0}) over threshold stopping times $\M^l_{\rm th}$. This implies (due to Theorem 2) (\ref{criteria110})--(\ref{criteria11}) and, besides, $p^*>I$. According to Theorem 3 $g(x)=x-I$ is $\r$-excessive function with respect to restriction of $X_t$ on the range $D^l_{p^*}$. This imply (\ref{criteria21}) --- see condition i) from the Statement.

Now, let (\ref{criteria110})--(\ref{criteria21}) hold. As in the proof of Theorem 3 consider the function
$$
\Phi(x)=V(p^*;x)=\left\{
\begin{array}{ll}
(p^*-I)\psi(x)/\psi(p^*), & \mbox{\rm for } x<p^*,\\[4pt]
x-I, &    \mbox{\rm for } x\ge p^*.
\end{array}
\right.
$$
Note, (\ref{criteria110}) imply $\Phi(x)\ge x-I$ for all $x\in \I$.

It is easy to check (similar to the corresponding arguments in the proof of Theorem 3) that all conditions i)--iii) of Lemma 2 are satisfied. Therefore, $\ds\Ex \Phi(X_{\t}) e^{-\r\t}\II_{\{\t{<}\infty\}}\le \Phi(x)$ for all $\t$ and $x\in \I$. Thus, $\ds\Ex (X_{\t}-I) e^{-\r\t}\II_{\{\t{<}\infty\}}\le \Phi(x)=V(p^*;x)$, i.e. $\t^l_{p^*}$ is optimal in the problem (\ref{optinv0}) over all stopping times.
\medskip

{\bf Remark.}
{\it As one can see from the proof of Theorem 5, conditions (\ref{criteria110})--(\ref{criteria21}) are necessary and sufficient for optimality of stopping time  $\t^l_{p^*}$ in optimal stopping problem with `optional' payoff function $g^+(x)=(x - I)^+$.}
\medskip

Let $X_t$ be geometric Brownian motion with rate of drift $\a$ and volatility $\s$. In this case, $\ds \psi (x)=x^{\b}$,\
where $\b$ is the positive root of the equation $ \frac12\ds \s^2 \b(\b-1)+\a\b=\r$. If $\r>\a$ then $\b>1$.
The optimal threshold is $\ds p^*=\frac{\b}{\b-1}I$ (see, e.g. \cite[Ch. 5]{DP}). Note, if $\s\to\infty$ then $\b\to 1$ and, therefore, $p^*\to\infty$. Hence, a large volatility increases threshold for investing.

\paragraph{\hspace*{\parindent} Optimal exit (abandonment) problem }

Another fundamental problem in Real options theory concerns finding the optimal time at which a decision-maker, who receives a payoff from an acting firm, should terminate operating and abandon firm.

Let $X_t$ be the price (at time $t$) of the good produced by the firm, and function $g_1(x)$ describes a dependence of firm's revenue on current price $x$.

The optimal abandonment problem is, using the available information about current production prices $X_t$, to find a moment $\t$ for terminating of firm's activity (exit moment) such that the net present value of the firm be maximal:
\beq\label{liq}
\E^x\left(\int_0^{\t}g_1(X_t) e^{-\r t}\,dt - L e^{-\r
\t}\right)\to \max_\t,
\eeq
where $L\ge 0$ is the abandonment cost.

The above Theorem 4 allows to rewrite the necessary and sufficient conditions for optimality of threshold stopping time $\t^r_{p^*}=\inf\{t{\ge} 0: X_t\le p^*\}$ (with assumption $L+R(x)\le 0$ for $x<p^*$) in the abandonment problem as follows:
\begin{eqnarray}
&&[L+R(p)]\varphi(p^*)\ge [L+R(p^*)]\varphi(p)\quad  \mbox{\rm for } p>p^*;\label{criteria201}\\[3pt]
&& I_2(p^*)S'(p^*)= L\varphi'(p^*); \label{criteria2011}\\[3pt]
&&  \L R(p) \ge \r R(p) +\r L \quad \mbox{\rm   for } p<p^*, \label{criteria22}
\end{eqnarray}
where functions $R(x)$ and $I_2(x)$ are specified in (\ref{Green}) and (\ref{I}) respectevely.

Let us note that $S'(p^*)>0,\, \varphi'(p^*)<0$. Therefore, if follows from (\ref{criteria2011}) that
\beq\label{neg}
I_2(p^*)=\int_{p^*}^r \varphi(y)g_1(y)H(y)\,dy<0.
\eeq

If we assume that function $g_1(x)$ is increasing, i.e. firm's revenue increases when current price rises, then from (\ref{neg}) we have immediately $g_1(p^*)<0$.  It means that a firm should terminate an activity only when its revenue will be strictly negative. In particular, if $g_1(x)=x-c$, where $x$ is current price of production, and $c$ is production cost, then optimal threshold of price $p^*$ for abandon a firm  should be less then production cost: $p^*<c$.

If a firm under consideration operates, e.g., in oil exploration, the optimal threshold $p^*$ can be interpreted as optimal (in sense of problem (\ref{liq})) breakeven oil price, that characterizes how low price can fall before oil projects start shutting down.

Consider the example (see \cite[Ch. 7]{DP}), where diffusion $X_t$ is geometric Brownian motion with rate of drift $\a<0$ and volatility $\s$, and $g_1(x)=x-c$. In this case,
$\ds
R(x)={x}/{(\r-\a)}-{c}/{\r},\quad \varphi (x)=x^{\b_1},\
$
where $\b_1$ is the negative root of the equation $ \frac12\ds \s^2 \b(\b-1)+\a\b=\r$.

If $c>\r L$, the optimal threshold is $\ds p^*=\frac{\b_1}{\b_1-1}\cdot\frac{\r-\a}{\r} (c-\r L)$. It is easy to see that $\b_1>\r/\a$, \  $p^*<c{-}\r L$ and all conditions (\ref{criteria201})--(\ref{criteria22}) hold. Note, that if volatility $\s$ tends to $+\infty$, then $\b_1\to 0$ and, therefore, $p^*\to 0$. It means that large volatility implies low level of price for shut down and more long period before abandonment of firm's activity, even though current profits $P_t-c$ will be negative.

\section{Free-boundary problem and\\ threshold strategies}

Let us return to a variational approach, described at Section 2, and give a new look to a free-boundary problem for the `threshold case'. Consider one-dimensional diffusion process $X_t$ with values in interval $D=(l,r)$ and the class of $l$-intervals. Assume that both drift and diffusion functions of $X_t$ are continuous, and payoff function is differentiable (on $D$).

For this case a free-boundary problem can be written as follows:
to find threshold $\bar p,\ l<\bar p<r$ \ and bounded twice differentiable function $H(x),\ l<x<\bar p$, such that
\begin{eqnarray}
&&\L H(x)=\r H(x), \quad l<x<\bar p;\label{fbp1} \\
&& H(\bar p-0)= g(\bar p), \label{fbp2}\\
&& H'(\bar p-0)=g'(\bar p)\label{fbp3}.
\end{eqnarray}

Conditions (\ref{fbp1})--(\ref{fbp2}) hold for the function
\beq\label{fbp0}
H(x)= h(\bar p)\psi(x),\quad   l<x<\bar p
\eeq
where $\psi(x)$ is an increasing solution to ODE (\ref{diffur}) and $h(p)=g(p)/\psi(p)$.
Thus, smooth pasting condition (\ref{fbp3}) at the point $\bar p$ \ is equivalent to $g(\bar p)\psi'(\bar p)=g'(\bar p)\psi(\bar p)$, or $h'(\bar p)=0$, i.e. $\bar p$ is a stationary point of the function $h$.

On the other hand, as the results of Section 3 show, the optimal threshold must be a point of maximum of the function $h$.
Of course, point of maximum is a stationary point of the same function, but not vice versa.

In economic literature (e.g. in Real Options theory) it is common opinion that solution to free-boundary problem is always a solution to optimal stopping problem.

To demonstrate a difference between solutions to optimal stopping problem and free-boundary problem, consider the following example.
\medskip

{\bfseries Example. \it A solution to free-boundary problem  can not give a solution to optimal stopping problem.}

\medskip

Let us consider geometric Brownian motion
$X_t=x\exp\{w_t\}$ (where $w_t$ be standard Wiener process),
payoff function $g(x){=} (x-1)^3+x^\d$ for $x\ge 0$, and discount rate $\r=\d^2/2$ ($\d>0$). The function $g$ is smooth and
increasing for all $\d>0$.

For this case the free-boundary problem is the following one:
\beq
\label{ex-Stefan}
\left\{
\begin{array}{l}
{\textstyle\frac12}x^2H''(x)+{\textstyle\frac12}xH'(x)=\r H(x),\qquad
0<x<p^* \\[6pt]
H(p^*)=g(p^*)\\[6pt]
H'(p^*)=g'(p^*)
\end{array}
\right.
\eeq

Remind the designations $V_p(x)=\Ex g(X_{\t_p})e^{-\r\t_p}$, and  $\t_p=\min\{t\ge 0: X_t\ge p\}$.

For $\d\le 3$ the free-boundary problem (\ref{ex-Stefan}) has the unique
solution: $ H(x)=V_1(x)=x^\d,\quad p^*=1$.  However, the stopping
time $\t_1$ (a first exit time over the level 1) is not optimal
since  $V_p(x)=((p-1)^3/p^\d+1)x^\d\to \infty$ when $p\to \infty$  (when
$\d<3$), and $V_p(x)\uparrow V(x)=2x^3$ (when $\d=3$) for any $x>0$.

For $\d> 3$ the free-boundary problem (\ref{ex-Stefan}) has two
solutions:
\begin{eqnarray}
&&  H(x){=}V_1(x){=} x^\d,\ p^*=1\quad \mbox{and} \label{Stef1}\\
&&
H(x){=}V_{p_\d}(x){=}h(p_\d)x^\d,\ p^*{=}p_\d{=}\d/(\d{-}3).\label{Stef2}
\end{eqnarray}

Note that $V_{p_\d}(x) > V_1(x)$. Thus, one of the solutions to free-boundary problem ((\ref{Stef1}),
corresponding to the boundary $p^*=1$) does not give a solution to
the optimal stopping problem (which there exists, in contrast to the
previous case).
\medskip

The above remarks together with classical results on extremal problems (especially, second-order optimality conditions) allow to derive new results about the relation between solutions to free-boundary problem (\ref{fbp1})--(\ref{fbp3}) and optimal stopping problem (\ref{optstop}) over the class $\M{=}\M^l_{\rm th}$ of threshold stopping times generated by all $l$-intervals.
\medskip

If function $g(x)$ be twice differentiable at the point $\bar p$, then representation (\ref{fbp0}) implies that $H''(\bar p{-}0)= h(\bar p)\psi''(\bar p)$. On the other hand, $g''(\bar p) = h''(\bar p)\psi(\bar p)+h(\bar p)\psi''(\bar p)$ since $h'(\bar p)=0$. Hence, we have
\beq
\label{2-order}
h''(\bar p)\psi(\bar p)=g''(\bar p)-H''(\bar p{-}0).
\eeq

Similarly, if  $h^{(k)}(\bar p)=0$, $k=1,...,n-1$ for some $n{>}2$ then
\beq
\label{n-order}
h^{(n)}(\bar p)\psi(\bar p)=g^{(n)}(\bar p)-H^{(n)}(\bar p{-}0).
\eeq

Due to Theorem 2 if $\t_{p^*}$ is optimal in the problem (\ref{optstop}) over the class $\M^l_{\rm th}$, then function $h(p)$ has maximum at the point $p^*$. Applying second-order condition for extremum and relation (\ref{2-order}) we obtain the following
\medskip

{\bf Proposition 1.} {\it
If $\t_{p^*}$ is optimal in the problem (\ref{optstop}) over the class $\M{=}\M^l_{\rm th}$ and function $g(x)$ is twice differentiable at the point $p^*$, then the
pair $(H(x),p^*)$, where $H(x)=h(p^*)\psi(x)$, is the solution to free-boundary problem (\ref{fbp1})--(\ref{fbp3}) and
 $H''(p^*{-}0)\ge g''(p^*)$.
}
\medskip

The inverse relation between solutions can be state as follows.
\medskip

{\bf Proposition 2.} {\it Let a pair $(H(x),p^*)$ be a solution to free-boundary problem (\ref{fbp1})--(\ref{fbp3}), such that $H''(p^*{-}0){>} g''(p^*)$ and one of the following conditions holds:

\emph{1)} $(H(x),p^*)$ is the unique solution to free-boundary problem (\ref{fbp1})--(\ref{fbp3});

\emph{2)}
if there exists another solution \ $(\tilde H(x),\tilde p)$   to free-boundary problem (\ref{fbp1})--(\ref{fbp3}), then $H(x)\ge g(x)$ for $l<x< p^*$, and either\\ \hspace*{1cm} \emph{(2a)} $\tilde p <p^*$; or\\  \hspace*{1cm} \emph{(2b)} $\tilde p >p^*$,
 $\tilde H^{(k)}(\tilde p{-}0){=}g^{(k)}(\tilde p)$, $k=2,...,n-1$, $\tilde H^{(n)}(\tilde p{-}0){>} g^{(n)}(\tilde p)$ for some $n{>}2$.

Then $\t_{p^*}$ is optimal in the problem (\ref{optstop}) over the class $\M{=}\M^l_{\rm th}$.
}
\medskip

{\sc Proof. } \
From the definition of $p^*$, inequality $H''(p^*{-}0){>} g''(p^*)$ and (\ref{2-order}) it follows that $h'(p^*)=0$ and $h'(p)$ strictly decreases at some neighborhood of $p^*$.

Let condition 1) holds.  It is easy to see that  $h'(p)>0$ for $p<p^*$ and $h'(p)<0$ for $p>p^*$ (else  $h'(q)=0$ for some $q\neq p^*$ and $(h(q)\psi(x), q)$ is another solution to free-boundary problem (\ref{fbp1})--(\ref{fbp3})). Then application of Theorem 2 gives the optimality of stopping time $\t_{p^*}$.

Now, let 2) holds. Then representation (\ref{fbp0})  immediately imply that $h(x)=g(x)/\psi(x)\le H(x)/\psi(x)=h(p^*)$ for $l<x< p^*$.

Let us prove that $h'(p)\le 0$ for all $p>p^*$. Indeed, if $h'(p_1)>0$ for some $p_1 >p^*$, then there exists $p^*<p_0<p_1$ such that $h'(p_0)=0$ and $h'(p)>0$ for all $p_0<p<p_1$.  Therefore, $(h(p_0)\psi(x), p_0)$ is another solution to free-boundary problem  (\ref{fbp1})--(\ref{fbp3}), and due to 2) $h^{(k)}(p_0)=0$, $k=2,...,n-1$, $h^{(n)}(p_0){<} 0$ for some $n{>}2$, that contradicts to positivity of $h'(p)$ for $p_0<p<p_1$.

Hence, $h'(p)\le 0$ for all $p>p^*$ and due to Theorem 2 $\t_{p^*}$ is optimal in the problem (\ref{optstop}) over the class $\M{=}\M^l_{\rm th}$ of $l$-threshold stopping times. \ $\square$
\medskip

Let us demonstrate, how Proposition 2 works for the example considered above.

In the case $\d\le 3$, we can't apply Proposition 2 since $H''(1){=}\d(\d-1){=}g''(1)$.

In the case $\d>3$, for the solution (\ref{Stef2}) with the boundary point $p^*=p_\d$ (see (\ref{Stef2})) we have
$$
H''(p_\d)=\d(\d-1)p_\d^{\d-2}+\frac{27(\d-1)}{\d(\d-3)}> \d(\d-1)p_\d^{\d-2}+\frac{18}{\d-3}=g''(p_\d).
$$
Besides,
$
\ds h'(x)=(x-1)^2x^{-\d-1}(\d-3)\left(\frac{\d}{\d-3}-x\right)\le 0
$
for $\d>3$ and $x<p_\d$, therefore, $H(x)=h(p_\d)x^\d\ge h(x)x^\d =g(x)$ for $x<p_\d$. Since for the another solution (\ref{Stef1}) boundary point   $p^*=1<p_\d$,  then the condition 2) of Proposition 2 holds.

Hence, $\t_{p_\d}$ is optimal in the problem (\ref{optstop}) over the class of threshold stopping times $\M^l_{\rm th}$.

Moreover, it can be shown (using Theorem 3), that  $\t_{p_\d}$ will be optimal stopping time over the class of all stopping times also. Indeed, by direct calculations one can obtain that
$$
\L g(x)-\r g(x)\le -1.5\frac{\d}{\d-3}<0 \qquad \mbox{\rm for } x>p_\d,
$$
and we can apply Theorem 3.

\end{document}